\documentclass[12pt]{article}
\usepackage{mathrsfs}
\usepackage{a4wide}
\usepackage[T1]{fontenc}
\usepackage{etoolbox}
\usepackage{amsmath,amsfonts,amssymb,amsthm}

\def \beq{\begin{equation}}
\def \eeq{\end{equation}}

\newcommand{\flux}{{\Phi^B}}
\newcommand{\fA}{{\phi^A}}
\def\X{\mathcal{X}}
\def\R{{\mathbb R}}
\def\beq{\begin{equation}}
\def\eeq{\end{equation}}

\def\f{\widetilde{f}_{(\delta)}}
\def\F{\widetilde{F}_{(\delta)}}
\def\bb1{{1\hspace*{-2.4pt}\rm{l}}}

\newtheorem{theorem}{Theorem}[section]
\newtheorem{definition}[theorem]{Definition}

\newtheorem{corollary}[theorem]{Corollary}

\theoremstyle{definition}
\newtheorem{remark}[theorem]{Remark}

\numberwithin{equation}{section}

\newcommand{\xx}{{\underline x}}

\newcommand{\y}{{\bf y}}

\newcommand{\x}{{\bf x}}

\newcommand{\la}{{\langle}}
\newcommand{\ra}{{\rangle}}

\begin{document}

\noindent 
\begin{center}
{\Large Spectral edge regularity of magnetic Hamiltonians}
\end{center}

\begin{center}
June 25, 2014
\end{center}

\begin{center}
 
Horia D. Cornean\footnote{Department of Mathematical Sciences,
  Aalborg University, Fredrik Bajers Vej 7G, 9220 Aalborg, Denmark} and Radu
Purice\footnote{Institute
of Mathematics Simion Stoilow of the Romanian Academy, P.O.  Box
1-764, Bucharest, RO-70700, Romania.}${{}^{,}}$\footnote{Laboratoire Europ\'een Associ\'e CNRS Franco-Roumain {\it Math-Mode}}

\end{center}

\vspace{0.2cm}

\begin{abstract}
We analyse the spectral edge regularity of a large class of magnetic Hamiltonians when the perturbation is generated by a globally bounded magnetic field. We can prove Lipschitz regularity of spectral edges if the magnetic field perturbation is either constant or slowly variable. We also recover an older result by G. Nenciu who proved Lipschitz regularity up to a logarithmic factor for general globally bounded magnetic field perturbations.
\end{abstract}

\section{Introduction and main results}

This is the second paper of the authors on the spectral regularity with respect to perturbations induced by Peierls-type magnetic flux phases. We assume that the flux is generated by a globally bounded magnetic field whose intensity is proportional with $\epsilon\in \R$. In a previous paper \cite{Co-Pu}, the regularity of the Hausdorff distance between the perturbed and unperturbed spectra was investigated.  In the current paper we analyse the regularity of spectral edges when $\epsilon$ varies.

It is well known that the magnetic perturbation induced by a non-decaying magnetic field is a singular perturbation and the spectral stability is not obvious. The first proof of spectral stability of
nearest-neighbor Harper operators with constant magnetic fields can be found in \cite{Ell} while in \cite{BEY} it is shown that the gap boundaries/spectral edges are
$\frac{1}{3}$-H\"older continuous in $\epsilon$. Later results \cite{AMS}, \cite{He-Sj1,
  He-Sj2}, \cite{HR} show that Hausdorff distance between spectra goes like
$|\epsilon-\epsilon_0|^{\frac{1}{2}}$. This result is optimal in the sense that it is known that gaps can appear/close down
precisely like $|\epsilon-\epsilon_0|^{\frac{1}{2}}$ if $\epsilon_0$ generates a rational flux or if the lattice is triangular, see \cite{He-Sj2, HKS, Bell1, BCS}.  

The first proof of Lipschitz continuity of gap edges 
for Harper-like operators with constant magnetic fields 
was given by Bellissard \cite{Bell2} (later on Kotani
\cite{Ko} extended his method to more general regular lattices and
dimensions larger than two).

In the continuous case of Schr\"odinger operators with bounded magnetic fields, the stability of gaps was first proved in \cite{AS} and \cite{Nenciu}. Then in \cite{BC} the H\"older exponent of gap edges was shown to be at least $\frac{2}{3}$, while \cite{C} provided a new proof of the results of \cite{Bell2} and extended them to continuous two dimensional Schr\"odinger operators perturbed by weak {\it constant} magnetic fields. We note that purely magnetic Schr\"odinger operators of Iwatsuka type (see \cite{DGR} and references therein) have magnetic bands whose width is proportional with the total variation of the magnetic field. An interesting open problem would be to see whether  such a behaviour remains true when the magnetic field is slightly perturbed around a non-zero constant value, and this perturbation is not a function of just one variable. 

A general discrete problem was formulated by Nenciu in \cite{Ne2}
where he worked with more general real and antisymmetric phases obeying a certain area condition (see \eqref{apr2}). These phases appear very naturally in the continuous case, see
\cite{CN, IMP1, LMS, MP1, MP2, Ne1}. 

Using a completely different method of proof, Nenciu showed in \cite{Ne2} that the gap edges are Lipschitz up to a
logarithmic factor. His method uses a regularity property of almost mid-convex
functions, and works for arbitrary bounded magnetic fields, not necessarily constant.

In the current paper we significantly improve our previous results in \cite{C, Co-Pu}. In particular, we recover the results of \cite{Ne2} and moreover, we can prove Lipschitz regularity of spectral edges if the magnetic field perturbation is either constant or slowly variable.  We also obtain results in the case in which the off-diagonal localization of the unperturbed kernels is weak. 

The structure of the paper is as follows. In the rest of this section we state the main results in Theorem \ref{theorem1} and Corollary \ref{corolar1}, where we also discuss which class of magnetic Hamiltonians/$\Psi$DO's are covered. The last two sections contain all the proofs.

\subsection{The main theorem}

We use the notation $\la \x \ra :=(1+|\x|^2)^{1/2}$ for any $\x\in\R^d$.
\begin{definition}\label{definition1}
 We say that a linear operator $T_{K}\in\mathbb{B}(L^2(\R^d))$ has an {\em  off-diagonal polynomial decay of order $\alpha\geq 0$}
if it is defined by an integral kernel $K$ such that:
 $$
  \|T_K\|_\alpha:=\max\left \{\underset{\x\in\R^d}{\sup}\int_{\R^d}\left|K(\x,\y)\right|\la \x-\y\ra ^\alpha\,d\y,\; \underset{\y\in\R^d}{\sup}\int_{\R^d}\left|K(\x,\y)\right|\la \x-\y\ra ^\alpha\,d\x\right \}<\infty.
 $$
We denote by $\mathscr{C}^\alpha$ the complex linear space of these operators with the norm $ \|\cdot \|_\alpha$, including by definition the identity operator ${\rm Id}$. 
\end{definition}
For us, a {\it magnetic field} will be described by a closed two-form on $\mathbb{R}^d$ of the form $B=dA$ with bounded components $B_{jk}(\x)=-B_{kj}(\x)=\partial_jA_k(\x)-\partial_kA_j(\x)\in\mathbb{\R}$. The quantity we shall be interested in, is the flux of this $2$-form through triangles (here $\la \x,\y,\x'\ra $ denotes the triangle with vertices  $\x$, $\y$ and $\x'$ in $\mathbb{R}^d$):
\begin{equation}\label{apr1}
\flux(\x,\y,\x'):= \int_{\la \x,\y,\x'\ra}B
=
\underset{j,k}{\sum}(y_j-x_j)(x'_k-y_k)\int_0^1dt\int_0^tds\,B_{jk}
(\x+t(\y-\x)+s(\x'-\y)).
\end{equation}
If we work with the {\it transverse gauge} which obeys $\sum\limits_{j=1}^{N}x_jA_j(\x)=0$, up to a use of the Stokes theorem we have that
\begin{equation}\label{apr1'}
\fA(\y,{\bf z}):=-\int_{[\y,{\bf z}]}A=-\flux(0,\y,{\bf z}) =-\underset{j<k}{\sum}(y_jz_k-y_kz_j)\int_0^1dt\int_0^tds\,B_{jk}
(t\y+s({\bf z}-\y)).
\end{equation}
In the case of a constant unit magnetic field we have $B_{jk}=-B_{kj}=1$ if $j<k$ and $B_{jj}=0$. Thus \eqref{apr1'} gives:
\begin{equation}\label{apr3}
\fA(\x,\y) =-\frac{1}{2}\sum_{j<k}(x_jy_k-x_ky_j).
\end{equation}

For every operator $T_K\in \mathscr{C}^\alpha \subset \mathbb{B}(L^2(\R^d))$ we can define a family of bounded linear operators $T_{K_A}\in\mathscr{C}^\alpha$ whose kernels are given by $$K_A(\x,\y):=e^{i\fA(\x,\y)}K(\x,\y).$$ 
More generally, one may consider a general $2$-point function $\phi:\R^d\times\R^d\rightarrow\mathbb{R}$ satisfying the following properties:
\begin{equation}\label{apr2}
\phi(\x,\y) =-\phi(\y,\x),\quad\text{and}\quad\left|\phi(\x,\y)+\phi(\y,{\bf z})+\phi({\bf z},\x)\right|\leq C |\la \x,\y, {\bf z}\ra|,
\end{equation}
where $C$ is some constant. It is easy to see that $\fA$ is antisymmetric and due to the Stokes theorem we have:
\begin{equation}\label{iunie26}
\fA(\x,\y)+\fA(\y,\x')-\fA(\x,\x')=-\flux(\x,\y,\x').
\end{equation}
Hence both conditions in \eqref{apr2} are satisfied if the components of $B$ obey $||B_{jk}||_\infty \leq C$. 

If one uses a different vector potential $A'$ such that $dA'=dA=B$, then Stokes theorem ensures that 
$$\fA(\x,{\bf y})= -\int_{[\x,{\bf y}]}A' +\int_{[0,{\bf y}]}A' -\int_{[0,\x]}A',$$
which shows that the operator with kernel 
$e^{-i\int_{[\x,{\bf y}]}A'}K(\x,\y)$ is unitarily equivalent with $T_{K_A}$. 

A slightly different but important situation is when the magnetic field perturbation comes from a slowly varying vector potential $A_\epsilon(\x):=A(\epsilon \x)$, where all the components of $dA$ are globally bounded together with their first order partial derivatives. In this case, the magnetic field perturbation is of the form $ \epsilon B_\epsilon(\x)$ with $B_\epsilon(\x):=(dA) (\epsilon \x)$. Then we define:
\begin{equation}\label{mai1}
\phi^{A_\epsilon}(\x,\x'):=-\epsilon \Phi^{B_\epsilon}(0,\x,\x'),\quad 
K_{A_\epsilon}(\x,\y):=e^{i\phi^{A_\epsilon}(\x,\y)}K(\x,\y).
\end{equation}

Let $\epsilon \in \R$ and $T_K\in\mathscr{C}^\alpha$ a self-adjoint operator. Since both $\phi^{\epsilon A}$ and $\phi^{A_\epsilon}$ are antisymmetric, both  $T_{K_{\epsilon A}}$ and $T_{K_{ A_\epsilon}}$ belong to $\mathscr{C}^\alpha$ and are self-adjoint. We are now ready to state our main theorem. 

\begin{theorem}\label{theorem1}
Let $T_K\in\mathscr{C}^\alpha$ be self-adjoint. Denote by $\mathcal{E}(\epsilon)=\sup\sigma(T_{K_{\epsilon A}})$ and by  $\mathcal{E}_\epsilon:=\sup\sigma(T_{K_{A_\epsilon }})$. The following statements hold true uniformly in $|\epsilon|\leq 1/2$: 
\begin{enumerate}
 \item If $1\leq\alpha<2$, then there exists a numerical constant $C_\alpha>0$ with $\lim_{\alpha\nearrow 2} C_\alpha=\infty$, such that
 \begin{equation}\label{gravi1}
 \left|\mathcal{E}(\epsilon)\,-\,\mathcal{E}(0)\right|\ \leq\ C_\alpha  \|T_K\|_\alpha \; |\epsilon|^{\alpha/2};
 \end{equation}
 \item If $\alpha\geq 2$, then there exists a numerical constant $C>0$  such that
 \begin{equation}\label{gravi2}
 \left|\mathcal{E}(\epsilon)\,-\,\mathcal{E}(0)\right|\ \leq\ C\|T_K\|_2 \; |\epsilon|\; \ln(1/|\epsilon|);
  \end{equation}
 \item Let $\alpha\geq 2$  and assume that either $B$ is a constant magnetic field, or the magnetic field perturbation comes from a slowly varying vector potential $A_\epsilon$.  Then there exists a numerical constant $C>0$ such that
 \begin{equation}\label{gravi3}
 \max\{ \left|\mathcal{E}(\epsilon)\,-\,\mathcal{E}(0)\right|, \left|\mathcal{E}_\epsilon-\mathcal{E}_0\right|\} \leq C\|T_K\|_2|\epsilon|;
  \end{equation}
\end{enumerate}
\end{theorem}

\vspace{0.2cm}

\begin{remark}\label{runu} Because $\inf(\sigma(T_{K_{ A}}))=-\sup(\sigma(T_{(-K)_{ A}}))$, the theorem also holds true if  $\mathcal{E}(\epsilon)=\inf\sigma(T_{K_{\epsilon A}})$ and $\mathcal{E}_\epsilon=\inf\sigma(T_{K_{A_\epsilon}})$. 
\end{remark}

\vspace{0.2cm}

\begin{remark}\label{rdoi}  We have the general identity 
$$||T_{K_{ A}}||=\max\{ -\inf(\sigma(T_{K_{ A}})), \sup(\sigma(T_{K_{ A}})) \}.$$
If $F(\epsilon):=\max\{f_1(\epsilon),f_2(\epsilon)\}$ let us show that $ |F(\epsilon)-F(0)|\leq \max_{n}|f_n(\epsilon)-f_n(0)|$. Assume without loss that $F(\epsilon)\geq F(0)$. If $F(\epsilon)=f_j(\epsilon)$ and $F(0)=f_j(0)$ then the inequality is trivial. If $F(\epsilon)=f_j(\epsilon)$ and $F(0)=f_k(0)$ with $j\neq k$ we have: 
$$|F(\epsilon)-F(0)|=F(\epsilon)-F(0)=f_j(\epsilon)-f_k(0)\leq f_j(\epsilon)-f_j(0)\leq \max_{n=1,2}|f_n(\epsilon)-f_n(0)|.$$
Thus the theorem also holds true if $\mathcal{E}(\epsilon)=||T_{K_{\epsilon A}}||$ and $\mathcal{E}_\epsilon =||T_{K_{A_\epsilon }}||$.
\end{remark}

\vspace{0.2cm}

\begin{remark}\label{rtrei} Since $\phi^{\epsilon A}=\epsilon\phi^A =(\epsilon-\epsilon_0)\phi^A +\epsilon_0\phi^A$, we can absorb $e^{i\epsilon_0 \phi^A(\x,\y)}$ into the kernel $K$ without changing its off-diagonal decay properties. Thus the above results for $\mathcal{E}(\epsilon)$ can be easily extended near any $\epsilon_0\in \R$ with $|\epsilon |$ replaced by $|\epsilon-\epsilon_0|$. The dependence on $\epsilon$ of $\phi^{A_\epsilon}$ is non-linear and it seems that the results on $\mathcal{E}_\epsilon$ cannot be extended. In fact, if $\epsilon$ is large, we can no longer talk about a slowly varying magnetic field. 
\end{remark}

\vspace{0.2cm}

\subsection{Application to magnetic pseudodifferential operators}

Let us briefly present the setting behind them. Denote by
$\mathcal{X}:=\mathbb{R}^d$ {\it the configuration space} of a
physical system, by $\X^*\cong\mathbb{R}^d$ its dual ({\it the space of
momenta}), by $\la \cdot,\cdot \ra:\X^*\times\X\rightarrow\R$ the
duality bilinear form and by $\Xi:=\X\times\X^*$ {\it the phase space} with the
canonical sympletic form
$\sigma\big((x,\xi),(y,\eta)\big):=\la \xi,y\ra -\la\eta,x\ra $. Let us recall from
\cite{MP1, MP2} that to any {\it classical Hamiltonian} described by a real
smooth function $h:\Xi\rightarrow\mathbb{R}$ (with polynomial growth together
with all its derivatives) and to any bounded smooth magnetic field described by
a closed 2-form having components $B_{jk}\in BC^\infty(\mathcal{X})$ we
associate a {\it quantum Hamiltonian} defined by the following action on
test functions (as oscillating integrals):
\beq
\big(\mathfrak{Op}^{A}(h)f\big)(\x):=(2\pi)^{
-d}\int_{\X}\int_{\X^*}e^{i\la \xi,\x-\y\ra}\Lambda^A(\x,\y)h\left (\frac{\x+\y}{2},
\xi\right )f(\y)d\y\,d\xi,
\eeq
where $A$ is a vector
potential such that $dA=B$, $\Lambda^A(\x,\y)\,:=e^{i\phi^{\epsilon A}(\x,\y)}$,  $f\in\mathscr{S}(\X)$, and $\x\in\X$. Let us
remind here that the quantum Hamiltonian
depends on the choice of the vector potential $A$, but different choices lead to
unitarily equivalent operators. Choosing a vector potential of class
$C^\infty_{\text{\sf pol}}(\X)$ is always possible and with such a choice,
Proposition 3.5 in \cite{MP1} states that the application $\mathfrak{Op}^A$
defines a bijection from the tempered distributions on $\Xi$ to the continuous
operators from $\mathscr{S}(\X)$ to $\mathscr{S}^\prime(\X)$. Thus the
composition of operators (when possible) induces a composition law on
$\mathscr{S}^\prime(\Xi)$, that we call {\it the magnetic Moyal product}, and is
 explicitly given by the following formula (for a pair of test functions $f$
and $g$ from $\mathscr{S}(\Xi)$) that only depends on the magnetic field and
not on the vector potential:
$$
\big(\phi\sharp^B\psi\big)(X)=\pi^{-2d}\int_\Xi\int_\Xi e^{
-2i\sigma(X-Y,X-Z)}e^{-i\int_{\mathcal{T}(\x,\y,{\bf z})}B}\phi(Y)\psi(Z)\,dY\,dZ
$$
where we have introduced the notation of the form
$X=(\x,\xi)$ for the points of $\Xi$ and we have denoted by $\mathcal{T}(\x,\y,{\bf z})$
the triangle with vertices $\x-\y+{\bf z}$, $\y-{\bf z}+\x$, and ${\bf z}-\x+\y$. 

Before stating the second main result of our paper, we need one more notation.  An interesting class of symbols, related to Onsager-Peierls effective
Hamiltonians, are symbols from $S^0_0(\Xi)$ that do not depend on the first
variable $x\in\X$ and are periodic in the second variable $\xi\in\X^*$ with
respect to a lattice $\Gamma_*\subset\X^*$. We denote by $S^0_{\Gamma_*}$
this class of symbols.

\begin{corollary}\label{corolar1}
Let $F$ be either a symbol in $S_1^t(\Xi)$ with $t<0$, or a symbol in $S^0_{\Gamma_*}$.  Let $\mathcal{E}(\epsilon)$ denote either $\sup \sigma (\mathfrak{Op}^{\epsilon A}(F))$, $\inf \sigma (\mathfrak{Op}^{\epsilon A}(F))$, or $||\mathfrak{Op}^{\epsilon A}(F)||$. Let $\mathcal{E}_\epsilon$ denote the same quantities defined with $A_\epsilon$ instead of $\epsilon A$. 

\noindent (i). There exists a constant $C<\infty$ such that:
$$|\mathcal{E}(\epsilon)-\mathcal{E}(0)|\leq C|\epsilon| \ln(1/|\epsilon|),\quad |\epsilon|\leq 1/2.$$

\noindent (ii). If the magnetic field is either constant or slowly variable, the logarithmic factor is absent:
$$\max\{|\mathcal{E}(\epsilon)-\mathcal{E}(0)|, |\mathcal{E}_\epsilon-\mathcal{E}_0|\}\leq C|\epsilon| ,\quad |\epsilon|\leq 1/2.$$
\end{corollary}

\vspace{0.2cm}

\begin{remark}\label{rpatru}
Let us recall from \cite{IMP1} that if $h$
is a real elliptic symbol (of H\"{o}rmander type), of strictly positive order
$m$, then the corresponding magnetic $\Psi$DO can be extended to a lower semibounded self-adjoint
operator denoted by $H_A$, acting in $L^2(\X)$ with domain {\it a
magnetic Sobolev space} (as defined also in \cite{IMP1}). Moreover, let us
recall that Proposition 6.5 from \cite{IMP2} implies the existence for any
$\mathfrak{z}\in\rho\big(H_A\big)$ of a symbol $r_B(h,\mathfrak{z})\in
S_1^{-m}(\X)$ such that
$\big(H_A-\mathfrak{z}\big)^{-1}=\mathfrak{Op}^{A}\big(r_B(h,\mathfrak{z}
)\big)$.
In
\cite{AMP} and \cite{Co-Pu} we proved that the spectrum of $H_{\epsilon A}$ varies continuously (as a subset of
$\mathbb{R}$) with the parameter $\epsilon$.
\end{remark}

\vspace{0.2cm}

\begin{remark}\label{rcinci}
If $h$ is as before and $\Phi\in C^\infty_0(\R)$, then by using the
Helffer-Sj\"{o}strand formula, it was proved in (\cite{IMP2}, Proposition 6.7) that
there exists a symbol $\Phi_B[h]\in S_1^{-m}(\X)$ such that the operator
$\Phi(H_A)$ defined by functional calculus with self-adjoint operators is in fact of the form $\Phi(H_A)=\mathfrak{Op}^A\big(\Phi_B[h]\big)$. 
The results in \cite{IMP2} imply that 
$$||\Phi(H_{\epsilon A})-\mathfrak{Op}^{\epsilon A}\big(\Phi_0[h]\big)||\leq C |\epsilon|,\quad \Phi_0[h]\in S_1^{-m}(\X),$$
which shows that the eventual non-Lipschitz behavior in the spectrum of $\Phi(H_{\epsilon A})$ can only come from the phase factor. We observe that $\mathfrak{Op}^{\epsilon A}\big(\Phi_0[h]\big)$ is covered by Corollary \ref{corolar1}. 
\end{remark}

\vspace{0.2cm}

\begin{remark}\label{rsase}
With $h$ as before, suppose that the Weyl quantized operator $H:=\mathfrak{Op}(h)$ has a bounded and isolated spectral island $\sigma_0$.  Then one can find a function
$\Phi\in C^\infty_0(\mathbb{R})$ such that $\Phi(t)=t$ on $\sigma_0$ and the support of $\Phi$ is disjoint from the rest of $\sigma(H)$.  Thus $\sigma(\Phi(H))=\sigma_0\cup \{0\}$. Up to a translation in energy we may suppose that $0\in \sigma_0$. It follows that $H_{\epsilon A}$ will still have an isolated spectral island $\sigma_\epsilon$ near $\sigma_0$ if $\epsilon$ is small enough, and its edges behave like in Corollary \ref{corolar1}. 
\end{remark}

\vspace{0.2cm}

\section{Proof of Theorem \ref{theorem1}}

Let us fix a non-zero and non-negative symmetric function $f\in C^\infty_0(\R^d)$. If $\delta>0$ and $\x\in\R^d$ we denote
by $f_\delta(\x):=f(\delta \x)$ and by
$\widetilde{f}_{(\delta)}:=f_\delta* f_\delta$. We have:
$$
0\leq \f(\x)\leq \|f_\delta\|_2^2=\int_{\R^d}d\y\,f_\delta(-\y)f_\delta(\y)=\f(0) =\delta^{-d}\|f\|_2^2.
$$

For any kernel $K$ such that $T_K\in \mathscr{C}^\alpha$ for some $\alpha\geq 0$ we define: 
\begin{equation}\label{aprilie14}
K_{(\delta)}(\x,\y):=K(\x,\y)\f(\x-\y)\|f_\delta\|_2^{-2}.
\end{equation}
Since $|K_{(\delta)}(\x,\y)|\leq |K(\x,\y)|$, we have 
$\left\|T_{K_{(\delta)}}\right\|\leq\left\| T_{K_{(\delta)}}\right\|_\alpha\leq \left\|T_K\right\|_\alpha$. 

Denote by $\F(\x,\y):=\f(\x-\y)$ and 
$F_{\delta}(\x,\y):=f_\delta(\x-\y)$ the kernels associated to the obvious  convolution operators. We now write two simple but very important identities. The first one is: 
\begin{align}\label{aprilie16'}
\|f_\delta\|_2^{-2}\int_{\R^d}|F_{\delta}(\y,\x)|^2d\y=\|f_\delta\|_2^{-2}\F(\x,\x)=\|f_\delta\|_2^{-2}\f(0)=1,\quad \forall \x\in\R^d.
\end{align}
The second one is (we use \eqref{iunie26}):
\begin{align}\label{aprilie16}
e^{i\phi^{\epsilon A}(\x,\x')}\F(\x,\x')&=
\int_{\R^d} e^{i\Phi^{\epsilon
B}(\x,\y,\x')}\left[e^{i\phi^{\epsilon
A}(\x,\y)}F_{\delta}(\x,\y)\right]\left[e^{i\phi^{\epsilon
A}(\y,\x')}F_{\delta}(\y,\x')\right]d\y\nonumber 
\\
&=\int_{\R^d}\,\left[\overline{e^{i\phi^{\epsilon
A}(\y,\x)}F_{\delta}(\y,\x)}\right]\left[e^{i\phi^{\epsilon
A}(\y,\x')}F_{\delta}(\y,\x')\right] d\y\nonumber 
\\
&+\int_{\R^d}\left[e^{i\Phi^{\epsilon
B}(\x,\y,\x')}-1\right]\left[e^{i\phi^{\epsilon
A}(\x,\y)}F_{\delta}(\x,\y)\right]\left[e^{i\phi^{\epsilon
A}(\y,\x')}F_{\delta}(\y,\x')\right]d\y\nonumber\\
&=\int_{\R^d}\left[\overline{e^{i\phi^{\epsilon
A}(\y,\x)}F_{\delta}(\y,\x)}\right]\left[e^{i\phi^{\epsilon
A}(\y,\x')}F_{\delta}(\y,\x')\right] d\y\nonumber 
\\
&+e^{i\phi^{\epsilon A}(\x,\x')}\int_{\R^d}\left[e^{i\Phi^{\epsilon
B}(\x,\y,\x')}-1\right]F_{\delta}(\x,\y)F_{\delta}(\y,\x')d\y\nonumber\\
&+e^{i\phi^{\epsilon A}(\x,\x')}\int_{\R^d}\left |e^{i\Phi^{\epsilon
B}(\x,\y,\x')}-1\right |^2F_{\delta}(\x,\y)F_{\delta}(\y,\x')d\y.
\end{align}
If we multiply the left hand side of \eqref{aprilie16'} with  $\|f_\delta\|_2^{-2}K(\x,\x')$ we obtain $e^{i\phi^{\epsilon A}(\x,\x')}K_{(\delta)}(\x,\x')$, see \eqref{aprilie14}. Then we can 
compute the quadratic form of $T_{K_{(\delta),\epsilon A}}$ on some $u\in L^2(\R^d)$ with
$\|u\|_{L^2(\R^d)}=1$:
\begin{align}\label{aprilie17}
&\langle u, T_{K_{(\delta),\epsilon A}}u\rangle=
\|f_\delta\|_2^{-2}\int_{\R^d}d\y\,\left\langle\left[e^{i\phi^{\epsilon
A}(\y,\cdot)}F_{\delta}(\y,\cdot)u(\cdot)\right]\,,\,T_K\left[e^{i\phi^{\epsilon
A}(\y,\cdot)}F_{\delta }(\y,\cdot)u(\cdot)\right]\right\rangle \nonumber \\
&+ \|f_\delta\|_2^{-2}\int_{\R^{d}}d\y \int_{\R^{d}}d\x \int_{\R^{d}}d\x'\left[e^{i\Phi^{
\epsilon
B}(\x,\y,\x')}-1\right]K_{\epsilon A}(\x,\x')F_{\delta}(\x,\y)F_{\delta}(\y,\x')\overline{u(\x)}u(\x')\nonumber \\
&+\|f_\delta\|_2^{-2}\int_{\R^{d}}d\y \int_{\R^{d}}d\x \int_{\R^{d}}d\x'\left |e^{i\Phi^{
\epsilon
B}(\x,\y,\x')}-1\right |^2K_{\epsilon A}(\x,\x')F_{\delta}(\x,\y)F_{\delta}(\y,\x')\overline{u(\x)}u(\x')\nonumber \\
&\leq  \mathcal{E}(0) \|f_\delta\|_2^{-2}\int_{\R^d}d\y\left\|e^{i\phi^{\epsilon
A}(\y,\cdot)}F_{\delta}(\y,\cdot)u(\cdot)\right\|^2_{L^2(\R^d)}\nonumber 
\\
&+ \|f_\delta\|_2^{-2}\int_{\R^{d}}d\y \int_{\R^{d}}d\x \int_{\R^{d}}d\x'\left[e^{i\Phi^{
\epsilon
B}(\x,\y,\x')}-1\right]K_{\epsilon A}(\x,\x')F_{\delta}(\x,\y)F_{\delta}(\y,\x')\overline{u(\x)}u(\x')\nonumber \\
&+\|f_\delta\|_2^{-2}\int_{\R^{d}}d\y \int_{\R^{d}}d\x \int_{\R^{d}}d\x'\left |e^{i\Phi^{
\epsilon
B}(\x,\y,\x')}-1\right |^2K_{\epsilon A}(\x,\x')F_{\delta}(\x,\y)F_{\delta}(\y,\x')\overline{u(\x)}u(\x').
\end{align}
Using \eqref{aprilie16'} in the last inequality we get: 
\begin{align}\label{aprilie17'}
&\langle u, T_{K_{(\delta),\pm \epsilon A}}u\rangle \leq  \mathcal{E}(0) \\
&+ \|f_\delta\|_2^{-2}\int_{\R^{d}}d\y \int_{\R^{d}}d\x \int_{\R^{d}}d\x'\left[e^{i\Phi^{
\pm \epsilon
B}(\x,\y,\x')}-1\right]K_{\pm \epsilon A}(\x,\x')F_{\delta}(\x,\y)F_{\delta}(\y,\x')\overline{u(\x)}u(\x')\nonumber \\
&+\|f_\delta\|_2^{-2}\int_{\R^{d}}d\y \int_{\R^{d}}d\x \int_{\R^{d}}d\x'\left |e^{i\Phi^{
\pm \epsilon
B}(\x,\y,\x')}-1\right |^2K_{\pm \epsilon A}(\x,\x')F_{\delta}(\x,\y)F_{\delta}(\y,\x')\overline{u(\x)}u(\x').\nonumber
\end{align}
Another simple but very important observation which we want to underline here, is that $T_K$ and $T_{K_{\pm \epsilon A}}$ have the same Schur $\alpha$-norms, see Definition \ref{definition1}. Moreover, 
we have the obvious identities:
\begin{align}\label{grav6}
K(\x,\x')=e^{i\phi^{\pm \epsilon A}(\x,\x')}
\left ( e^{ i\phi^{\mp \epsilon A}(\x,\x')}K(\x,\x')\right ),\quad 
T_K=T_{(K_{\mp \epsilon A})_{\pm \epsilon A}}. 
\end{align}
Thus changing $K$ with $K_{\mp \epsilon A}$ in \eqref{aprilie17'} we obtain:
\begin{align}\label{aprilie17''}
&\langle u, T_{K_{(\delta)}}u\rangle \leq  \mathcal{E}(\mp \epsilon )  \\
&+ \|f_\delta\|_2^{-2}\int_{\R^{d}}d\y \int_{\R^{d}}d\x \int_{\R^{d}}d\x'\left[e^{i\Phi^{
\pm \epsilon
B}(\x,\y,\x')}-1\right]K(\x,\x')F_{\delta}(\x,\y)F_{\delta}(\y,\x')\overline{u(\x)}u(\x')\nonumber \\
&+4\|f_\delta\|_2^{-2}\int_{\R^{d}}d\y \int_{\R^{d}}d\x \int_{\R^{d}}d\x'\left |\sin(\Phi^{|\epsilon|
B}(\x,\y,\x')/2)\right |^2K(\x,\x')F_{\delta}(\x,\y)F_{\delta}(\y,\x')\overline{u(\x)}u(\x').\nonumber 
\end{align}
This implies: 
\begin{align}\label{aprilie17'''}
&\langle u, T_{K_{(\delta)}}u\rangle \leq  \frac{1}{2}\mathcal{E}(\epsilon ) + \frac{1}{2}\mathcal{E}(-\epsilon )\\
&+ \|f_\delta\|_2^{-2}\int_{\R^{d}}d\y \int_{\R^{d}}d\x \int_{\R^{d}}d\x'\left (\cos(\Phi^{|\epsilon|
B}(\x,\y,\x'))-1\right )
K(\x,\x')F_{\delta}(\x,\y)F_{\delta}(\y,\x')\overline{u(\x)}u(\x')\nonumber \\
&+4\|f_\delta\|_2^{-2}\int_{\R^{d}}d\y \int_{\R^{d}}d\x \int_{\R^{d}}d\x'\left |\sin(\Phi^{|\epsilon|
B}(\x,\y,\x')/2)\right |^2K(\x,\x')F_{\delta}(\x,\y)F_{\delta}(\y,\x')\overline{u(\x)}u(\x').\nonumber 
\end{align}

Taking the supremum with respect to $u\in L^2(\R^d)$ with $\|u\|_{L^2(\R^d)}=1$
we obtain
\begin{align}\label{aprilie18}
&\sup (\sigma (T_{K_{(\delta),\epsilon A}}))\leq \frac{1}{2}\mathcal{E}(\epsilon ) + \frac{1}{2}\mathcal{E}(-\epsilon )+ \|f_\delta\|_2^{-2}\\
&\times \left (\underset{\|u\|=1}{\sup}\int_{\R^{3d}}d\y d\x d\x'\left |\cos(\Phi^{|\epsilon|
B}(\x,\y,\x'))-1\right |
 |K(\x,\x')| F_{\delta}(\x,\y)F_{\delta}(\y,\x')|u(\x)||u(\x')|\right . \nonumber \\
&\left .+4\underset{\|u\|=1}{\sup}\int_{\R^{3d}}d\y d\x d\x'\left |\sin(\Phi^{|\epsilon|
B}(\x,\y,\x')/2)\right |^2|K(\x,\x')|F_{\delta}(\x,\y)F_{\delta}(\y,\x')|u(\x)||u(\x')|\right ).\nonumber 
\end{align}
Using \eqref{apr1} and \eqref{apr2}, we obtain that for all $\x,\y,\x'\in \R^d$:
\begin{equation}\label{aprilie19}
|\Phi^{\epsilon B}(\x,\y,\x')|\leq C|\epsilon| |\langle \x,\y,\x'\ra|\leq \frac{C|\epsilon|}{2}\; |\x-\x'|\; |\x-\y|^{1/2}\; |\y-\x'|^{1/2}.
\end{equation}
From now on, $C$ will denote a generic numerical constant and not just the upper bound on the magnetic field. Using \eqref{aprilie19} and elementary properties of $\sin$ and $\cos$, we can bound the last two terms from the right hand side of \eqref{aprilie18} by:
\begin{align*}
&C  \|f_\delta\|_2^{-2}\underset{\|u\|=1}{\sup}\int_{\R^{3d}}|\Phi^{\epsilon B}(\x,\y,\x')|^{\alpha} 
|K(\x,\x')|F_{\delta}(\x,\y)F_{\delta}(\y,\x')|u(\x)|\; |u(\x')|d\y d\x d\x'\\
&\leq C |\epsilon|^\alpha \|T_K\|_\alpha\|f_\delta\|_2^{-2}\left(\int_{\R^d}
|\y|^{\alpha}f_\delta(\y)^2d\y\right) \leq C|\epsilon|^\alpha \delta^{-\alpha}\|T_K\|_\alpha ,\quad 1\leq \alpha <2,
\end{align*}
where the last inequality is due to the fact that $|\y|\leq C\delta^{-1}$ on the support of $f_\delta$. Moreover, 
\begin{align*}
&C  \|f_\delta\|_2^{-2}\underset{\|u\|=1}{\sup}\int_{\R^{3d}}|\Phi^{\epsilon B}(\x,\y,\x')|^{2} 
|K(\x,\x')|F_{\delta}(\x,\y)F_{\delta}(\y,\x')|u(\x)|\; |u(\x')|d\y d\x d\x'\\
&\leq C |\epsilon|^2 \|T_K\|_2\|f_\delta\|_2^{-2}\left(\int_{\R^d}
|\y|^{2}f_\delta(\y)^2d\y\right) \leq C|\epsilon|^2 \delta^{-2}\|T_K\|_2 ,\quad \alpha \geq 2.
\end{align*}
 Using this in \eqref{aprilie18} we obtain:
\begin{align}\label{aprilie20}
\sup (\sigma (T_{K_{(\delta)}}))&\leq \frac{1}{2}\mathcal{E}(\epsilon ) + \frac{1}{2}\mathcal{E}(-\epsilon )
+C|\epsilon|^\alpha\delta^{-\alpha}\|T_K\|_\alpha,\quad 1\leq \alpha<2, \nonumber \\
\sup (\sigma (T_{K_{(\delta)}}))&\leq \frac{1}{2}\mathcal{E}(\epsilon ) + \frac{1}{2}\mathcal{E}(-\epsilon )
+C|\epsilon|^2\delta^{-2}\|T_K\|_2,\quad \alpha\geq 2.
\end{align}
In general, if $A$ and $B$ are two bounded and self-adjoint operators and $u$ a norm-one vector we have:
$$\la u,A u\ra\leq \la u, B u\ra\ +||A-B||\leq \sup (\sigma (B))+||A-B||$$
which leads to:
$$ \sup (\sigma (A))\leq \sup (\sigma (B))+||A-B||.$$
Applying this inequality with $A=T_{K}$ and $B=T_{K_{(\delta)}}$, and using \eqref{aprilie20} we have:
\begin{align}\label{aprilie21}
\mathcal{E}(0)&\leq \frac{1}{2}\mathcal{E}(\epsilon ) + \frac{1}{2}\mathcal{E}(-\epsilon )
+C|\epsilon|^\alpha\delta^{-\alpha}\|T_K\|_\alpha+\left\|T_{K_{(\delta)}}-T_{K}\right\|,\quad 1\leq \alpha<2, \nonumber \\
\mathcal{E}(0)&\leq \frac{1}{2}\mathcal{E}(\epsilon ) + \frac{1}{2}\mathcal{E}(-\epsilon )
+C|\epsilon|^2\delta^{-2}\|T_K\|_2+\left\|T_{K_{(\delta)}}-T_{K}\right\|,\quad \alpha\geq 2.
\end{align}
Because $T_{K_{(\delta)}}-T_{K}=T_{K_{(\delta)}-K}$, using the Schur-Holmgren bound we have:
\begin{align}\label{grav3}
\left\|T_{K_{(\delta)}}-T_{K}\right\|&\leq \|f_\delta\|_2^{-2}
\underset{\x\in\R^d}{\sup}\int_{\R^d}
\left |\f(\x-\x')-\f(0)\right |  |K(\x,\x')|d\x'
\end{align}
The following bound is valid for all $\x\in\R^d$:
\begin{align*}
|\f(\x)-\f(0)|\leq \delta |\x| \int_0^1\int_{\R^d}|\nabla f(\delta t\x-\delta \y)|\; |f(\delta \y)| d\y dt\leq C\; |\x|\; \delta^{1-d}.
\end{align*}
Moreover, because $\f$ has a maximum at $\x=0$, by expanding up to the second order around $\x=0$ we have:
\begin{align*}
|\f(\x)-\f(0)|\leq \;\delta^2 |\x|^2 \int_0^1dt \int_0^tds \int_{\R^d} \sqrt{\sum_{1\leq j,k\leq d}| \partial_{jk}^2f(\delta s\x-\delta \y)|^2} \; |f(\delta \y)| d\y
\end{align*}
which leads to:
\begin{align}\label{grav2}
|\f(\x)-\f(0)|\leq C\; |\x|^2  \; \delta^{2-d}.
\end{align}
If $1\leq \alpha\leq 2$ we can combine the last two inequalities and get: 
\begin{align}\label{grav1}
|\f(\x)-\f(0)|\leq C\; |\x|^{\alpha}  \; \delta^{\alpha-d}.
\end{align}
Introducing \eqref{grav1} in \eqref{grav3} we obtain:
\begin{align}\label{grav4}
\left\|T_{K_{(\delta)}}-T_{K}\right\|&\leq C\; \delta^d 
\underset{\x\in\R^d}{\sup}\int_{\R^d}
|\x-\x'|^\alpha  \delta^{\alpha-d} |K(\x,\x')|d\x'\leq  C\delta^\alpha ||T_K||_\alpha ,\quad 1\leq \alpha\leq 2.
\end{align}
If $T_K\in \mathscr{C}_\alpha$ with $\alpha\geq 2$, we introduce \eqref{grav2} into \eqref{grav3} and obtain: 
\begin{align}\label{grav5}
\left\|T_{K_{(\delta)}}-T_{K}\right\|&\leq   C\delta^2 ||T_K||_2,\quad \alpha\geq 2.
\end{align}
Introducing the last two estimates in \eqref{aprilie21} we obtain
$$
\mathcal{E}(0)\leq \frac{1}{2}\mathcal{E}(\epsilon ) + \frac{1}{2}\mathcal{E}(-\epsilon )+C\|T_K\|_{\alpha}\left(\delta^\alpha
+|\epsilon|^\alpha\delta^{-\alpha}\right),\quad 1\leq \alpha < 2,
$$ 
or 
$$
\mathcal{E}(0)\leq \frac{1}{2}\mathcal{E}(\epsilon ) + \frac{1}{2}\mathcal{E}(-\epsilon )+C\|T_K\|_{2}\left(\delta^2
+|\epsilon|^2\delta^{-2}\right),\quad \alpha \geq 2.
$$ 
Until now, $\delta$ and $\epsilon$ have been independent. Keeping $\epsilon$ fixed and minimizing the right hand side with respect to $\delta$ imposes the condition 
$\delta= |\epsilon| ^{1/2}$ in both cases. Thus we obtain: 
\begin{align}\label{grav5'}
&\mathcal{E}(0)\leq\frac{1}{2}\mathcal{E}(\epsilon ) + \frac{1}{2}\mathcal{E}(-\epsilon )+C\|T_K\|_{\alpha}|\epsilon|^{\alpha/2},\quad 1\leq \alpha< 2,\nonumber \\
&\mathcal{E}(0)\leq\frac{1}{2}\mathcal{E}(\epsilon ) + \frac{1}{2}\mathcal{E}(-\epsilon )+C\|T_K\|_{2}|\epsilon|,\quad \alpha\geq 2. 
\end{align}
As we commented in Remark \ref{rtrei}, the above estimates can be obtained near any $\epsilon_0$ by redefining $K$. Thus we have just proved that the map $\R\ni x\mapsto \mathcal{E}(x)\in \R$ is a bounded, almost mid-convex function which obeys:
\begin{align}\label{grav7}
&\mathcal{E}\left (\frac{a+b}{2}\right )\leq\frac{1}{2}\mathcal{E}(a) + \frac{1}{2}\mathcal{E}(b)+C \|T_K\|_{2\beta} \left |\frac{b-a}{2}\right |^{\beta},\quad 1/2\leq \beta:=\alpha/2< 1, \\
&\mathcal{E}\left (\frac{a+b}{2}\right )\leq\frac{1}{2}\mathcal{E}(a) + \frac{1}{2}\mathcal{E}(b)+C\|T_K\|_{2}\left |\frac{b-a}{2}\right |,\quad \alpha\geq 2. \nonumber 
\end{align}

\vspace{0.2cm}

\noindent {\bf Regularity of bounded and almost mid-convex functions.} Now we shall prove that \eqref{grav7} implies \eqref{gravi1}, essentially following Nenciu \cite{Ne2}. Assume that $1/2\leq \beta <1$ is fixed and denote by $M:= C \|T_K\|_{2\beta}$. Thus we have:
$$\mathcal{E}\left (\frac{a+b}{2}\right )\leq\frac{1}{2}\mathcal{E}(a) + \frac{1}{2}\mathcal{E}(b)+M \left |\frac{b-a}{2}\right |^{\beta},\quad \forall a,b\in\R.$$
The strategy is to construct a constant $C_\beta>0$ such that for every $x\in \R$ and $0<\eta<1/2$ to have:  
\begin{align}\label{grav8}
-C_\beta \eta^\beta\leq \mathcal{E}(x+\eta)-\mathcal{E}(x)\leq C_\beta \eta^\beta.
\end{align}
One can easily prove by induction the following two inequalities (assume that $a<b$ and $n\geq 1$):
 \begin{align}\label{grav9}
\mathcal{E}\left ((2^{-1}+\dots +2^{-n})a+ 2^{-n}b\right )&\leq 
(2^{-1}+\dots +2^{-n})\mathcal{E}(a)+2^{-n}\mathcal{E}(b) \nonumber \\
&+ M\left (\frac{b-a}{2^n}\right )^\beta\left (1+\frac{1}{2^{1-\beta}}+\dots +\frac{1}{(2^{1-\beta})^{n-1}}\right )
\end{align}
and 
\begin{align}\label{grav10}
\mathcal{E}\left (2^{-n}a+ (2^{-1}+\dots +2^{-n})b\right )&\leq 
2^{-n}\mathcal{E}(a)+(2^{-1}+\dots +2^{-n})\mathcal{E}(b)
\nonumber \\
&+ M\left (\frac{b-a}{2^n}\right )^\beta\left (1+\frac{1}{2^{1-\beta}}+\dots +\frac{1}{(2^{1-\beta})^{n-1}}\right )
\end{align}
Given $\eta\in (0,1/2)$ we define $N_\eta:=\left \lfloor\frac{\ln(1/\eta)}{\ln(2)}\right \rfloor +1$. We have $N_\eta -1\leq \frac{\ln(1/\eta)}{\ln(2)}<N_\eta$, i.e. $1< \eta 2^{N_\eta}\leq 2$. 
Replace $a=x$, $b=x+\eta 2^{N_\eta}$ and $n=N_\eta$ in \eqref{grav9}. We have $(2^{-1}+\dots +2^{-N_\eta})a+ 2^{-N_\eta}b=x+\eta$ and:
$$\mathcal{E}(x+\eta)-\mathcal{E}(x)\leq \eta \frac{\mathcal{E}(b)-\mathcal{E}(a)}{b-a}+M\eta^\beta \frac{1}{1-2^{\beta -1}}.$$
Since $\mathcal{E}$ is bounded and we always have  $1<b-a\leq 2$, the right hand side is of order $\eta^\beta$.  Thus the right hand side of \eqref{grav8} is proved. 

For the other inequality in \eqref{grav8}, we replace $a=x+\eta -\eta 2^{N_\eta}$, $b=x+\eta$ and $n=N_\eta$ in \eqref{grav10}. We have $2^{-N_\eta}a+(2^{-1}+\dots +2^{-N_\eta})b=x$ and:
$$\mathcal{E}(x)-\mathcal{E}(x+\eta)\leq -\eta \frac{\mathcal{E}(b)-\mathcal{E}(a)}{b-a}+M\eta^\beta \frac{1}{1-2^{\beta -1}}.$$
Now we have to change sign and note again that $1<b-a\leq 2$, uniformly in $\eta$. Thus \eqref{grav8} is proved, and so is the first part of our Theorem \ref{theorem1}.  

Concerning the second part, i.e. the estimate \eqref{gravi2}, we see that 
both \eqref{grav9} and \eqref{grav10} hold true even if $\beta=1$. In this case, we can no longer use the geometric series and we get an extra $N_\eta$. This is the reason for having the logarithmic factor in \eqref{gravi2}. We give no further details. 

\vspace{0.2cm}

\noindent{\bf Constant magnetic field.} Now let us separately treat the case in which the perturbation comes from a constant magnetic field and $\alpha\geq 2$. In this case we shall see that one can directly prove a Lipschitz regularity for $\mathcal{E}$, without the logarithmic factor, and without using the trick based on almost mid-convex functions. 

Going back to the inequality \eqref{aprilie17'} we see that we can isolate the $\y$ integral in the second term on the right hand side. This integral is:
$$\int_{\R^d} F_\delta(\x,\y)F_\delta(\y,\x')(e^{i\Phi^{\pm \epsilon B}(\x,\y,\x')}-1)d\y.$$
Now let us show that the first order term in $\epsilon$ is just zero:
 $$\int_{\R^d} F_\delta(\x,\y)F_\delta(\y,\x')\Phi^{\pm \epsilon B}(\x,\y,\x')d\y=0.$$
Let us start by noticing that the above integral is proportional with:
 $$\sum_{j,k}B_{jk}\int_{\R^d} (x_k-y_k)f_\delta(\x-\y)(y_j-x'_j)f_\delta(\y'-\x')d\y. $$
Denote by $\mathcal{F}$ and $\mathcal{F}^-$ the Fourier transform and its
inverse. Then $$(x_k-y_k)f_\delta(\x-\y)\sim [\mathcal{F}^-(\partial_k \mathcal{F}f_\delta)](\x-\y),\quad (y_j-x_j')f_\delta(\y-\x')\sim [\mathcal{F}^-(\partial_j \mathcal{F}f_\delta)](\y-\x').$$
Hence:
\begin{align*}
\int_{\R^d} (x_k-y_k)f_\delta(\x-\y)(y_j-x'_j)f_\delta(\y'-\x')d\y &\sim  \{[\mathcal{F}^-(\partial_k \mathcal{F}f_\delta)]*[\mathcal{F}^-( \partial_j \mathcal{F}f_\delta)]\}(\x-\x')\\
&\sim \{ \mathcal{F}^-[(\partial_k \mathcal{F}f_\delta)(\partial_j \mathcal{F}f_\delta))]\}(\x-\x').
\end{align*}
But the product $(\partial_k \mathcal{F}f_\delta)(\partial_j \mathcal{F}f_\delta))$ is symmetric in $k$ and $j$ while $B_{jk}$ is antisymmetric, hence the sum gives zero.

Thus we see that in this case, the linear term disappears without having to appeal to the arithmetic mean, as we did in \eqref{aprilie18}. By performing the same type of analysis as before in order to deal with the quadratic terms, we obtain in particular that 
$$ |\mathcal{E}(\epsilon)-\mathcal{E}(0)|\leq C ||T||_2(\delta^2+|\epsilon|^2/\delta^2).$$
Choose $\delta =|\epsilon|^{1/2}$ and the proof is over. 

\vspace{0.2cm}

\noindent{\bf Slowly varying magnetic field.} In this case, the antisymmetric form entering the flux formula \eqref{apr1} is of the form $ B(\epsilon \; \cdot )$, while the total magnetic field perturbation is $ \epsilon B(\epsilon \; \cdot)$. Also, $\phi^{\pm\epsilon A}$ has to be replaced with $\phi^{\pm A_{\epsilon}}$. Again, the only obstacle in getting the Lipschitz behavior is the linear term as before. Let us show that we have the bound: 
\begin{equation}\label{gravi11}
\left | \int_{\R^d} F_\delta(\x,\y)F_\delta(\y,\x')\Phi^{\pm \epsilon B_\epsilon }(\x,\y,\x')d\y\right |\leq C |\x-\x'|\; |\epsilon|^2/\delta^2.
\end{equation}
Indeed, using \eqref{apr1} and Taylor's formula we have: 
$$\left |\Phi^{\epsilon B_\epsilon}(\x,\y,\x')-\frac{\epsilon}{2} \sum_{j,k}B_{jk}(\epsilon \x)(y_j-x_j)(x_k'-y_k)\right |\leq C |\epsilon|^2 |\la \x,\y,\x'\ra | (|\y-\x|+|\y-\x'|).$$
The contribution coming from $\epsilon \sum_{j,k}B_{jk}(\epsilon \x)(y_j-x_j)(x_k'-y_k)$ is zero, as in the constant case. The right hand side can be bounded by:
$$C |\epsilon|^2 |\x-\x'| |\x-\y|(|\y-\x|+|\y-\x'|),$$
term having a polynomial growth which introduced in the integral will generate a  diverging factor $\delta^{-2}$. Note that $|\x-\x'|$ can be coupled later on with $K(\x,\x')$. Having proved \eqref{gravi11}, the estimate we get in the end is: 
$$ |\mathcal{E}_\epsilon-\mathcal{E}_0|\leq C ||T||_1 |\epsilon|^2/\delta^2 + C ||T||_2(\delta^2+|\epsilon|^2/\delta^2)$$
 which gives the Lipschitz regularity by again taking $\delta=|\epsilon|^{1/2}$.
 \qed

 \section{Proof of Corollary \ref{corolar1}}
 
It was proved in 
\cite{MP1} that the magnetic quantization associated to the vector potential
$A$ is a topological vector space isomorphism
$\mathscr{S}^\prime(\Xi)\rightarrow\mathbb{B}\big(\mathscr{S}(\X);\mathscr{S}
^\prime(\X^*)\big)$. We have also given the explicit form of this isomorphism
by constructing the distribution kernel associated to a symbol. More precisely,
let us denote by $S_W:\X^2\rightarrow\X^2$ the linear isomorphism
$S_W(\x,\y):=\left(\frac{\x+\y}{2},\x-\y\right)$, by
$S_W^*:\mathscr{S}^\prime(\X^2)\rightarrow\mathscr{S}^\prime(\X^2)$ its
transposed map $S_W^*(F):=F\circ S_W$ and by
$\mathcal{F}:\mathscr{S}^\prime(\X)\rightarrow\mathscr{S}^\prime(\X^*)$ the
Fourier transform (normed in order to give a unitary map $L^2(\X)\rightarrow
L^2(\X^*
)$); we shall denote its inverse by $\mathcal{F}^-$. Then the map
$\mathfrak{K}_W:=S_W^*\circ(\bb1\otimes\mathcal{F}^-):\mathscr{S}
^\prime(\Xi)\rightarrow\mathscr{S}^\prime(\X^2)$ is a bijection associating to
any "symbol" on $\Xi$ an "integral kernel" on $\X$. 

We denote by
$T_K:\mathscr{S}(\X)\rightarrow\mathscr{S}^\prime(\X)$ the operator
associated to the integral kernel $K\in\mathscr{S}^\prime(\X^2)$, i.e.:
$$
\left\langle v, T_Ku\right\rangle\ :=\ K(v\otimes
u),\qquad\forall v,u \in \mathscr{S}(\X),
$$
or formally:
$$
\big(T_Ku\big)(\x)\ :=\ \int_\X K(\x,{\bf z})u({\bf z})d{\bf z}.
$$
Then we have the equality:
$$
\mathfrak{Op}(F)\ =\
T_{S_W^*\circ(\bb1\otimes\mathcal{F}^-)F},\quad\forall
F\in\mathscr{S}^\prime(\Xi).
$$
We make the important observation that the magnetic quantization can be expressed as:
$$
\mathfrak{Op}^A(F) =
T_{e^{i\phi^A}S_W^*\circ(\bb1\otimes\mathcal{F}^-)F},
$$
where $\phi^A(\x,\y)=-\int_{[\x,\y]}A$. If
$K:=S_W^*\circ(\bb1\otimes\mathcal{F}^-)F$, then the "magnetic integral
kernel" of \cite{Ne1, Ne2, C} is:
$$
K_A(\x,\y) = e^{i\phi^A(\x,\y)}K(\x,\y).
$$
Thus we have an explicit way of transferring results and formulas between
the two representations, working with the one which is more suitable for a given problem. The magnetic pseudodifferential calculus developed
in
\cite{MP1,IMP1,IMP2,MPR} is an equivalent formulation of the
calculus with magnetic integral kernels proposed in \cite{Ne1,Ne2, C}, the
equivalence being realized through the application taking a symbol into the
distribution kernel associated to the pseudodifferential operator the symbol generates.

\subsection{Decaying symbols}

Using Proposition 1.3.3 from \cite{ABG} and its variant given
in \cite{MPR}, we see that for any symbol $F$ of the type $S^{t}_1(\Xi)$ with
$t<0$, its partial inverse Fourier transform $(\bb1\otimes\mathcal{F}^-)F$
is a function for which there exists a constant $C$ such that
$$\sup_{\x\in \R^d}|[(\bb1\otimes\mathcal{F}^-)F](\x,\x')| \leq C |\x'|^{-d-t},\quad \x'\neq 0,$$
and has rapid decrease in the second variable (thus in $\x-\y$ for the kernel). Thus
through our identification
discussed above, it defines an integral operator with a kernel of class $\mathscr{C}^N$ for any $N\in\mathbb{N}$ (see Definition \ref{definition1}). Thus for this class of symbols, the Corollary is an immediate consequence of Theorem \ref{theorem1}.

\subsection{Periodic symbols}

 For any $\lambda\in S^0_{\Gamma_*}$ we denote by
$\tilde{\lambda}:=\big(\bb1\otimes\mathcal{F}^-\big)\lambda$ and taking into
account the Theorem in \cite{Ho-4} concerning the Fourier transform of
periodic distributions and denoting by $\Gamma\subset\X$ the dual lattice of
$\Gamma_*$, we obtain
\beq
\big(\mathfrak{Op}^A(\lambda)u\big)(\x)=\underset{\gamma''\in\Gamma}{\sum}
\Lambda^A(\x,\x-\gamma'')\tilde{\lambda}(\gamma'')u(\x-\gamma'')
\eeq
and the operator $\mathfrak{Op}^A(\lambda)$ has the distribution kernel
\beq
K^A_\lambda(\x,\y):=\underset{\gamma''\in\Gamma}{\sum}e^{i\phi^A(\x,\y)}\tilde{
\lambda} (\gamma'')\delta(\x-\y-\gamma'')
\eeq
with $\tilde{\lambda}(\gamma)$ having rapid decay with respect to
$\gamma\in \Gamma$. 

There exists a $d$-dimensional parallelepiped $\Omega$ such that every $\x\in \R^d$ can be uniquely represented as $\gamma +\xx$, with $\gamma\in \Gamma$ and $\xx\in \Omega$. We can see $\mathfrak{Op}^{\epsilon A}(\lambda)$ as an operator in $$l^2(\Gamma;L^2(\Omega))\sim l^2(\Gamma)\otimes L^2(\Omega),\quad \mathfrak{Op}^A(\lambda) =\{T_{\gamma \gamma'}\}_{\gamma,\gamma'\in \Gamma},\quad T_{\gamma \gamma'}\in \mathbb{B}(L^2(\Omega)),$$
where the operator $T_{\gamma \gamma'}$ has the distribution kernel:
$$T_{\gamma \gamma'}(\xx,\xx'):=K^{\epsilon A}_\lambda(\gamma +\xx,\gamma' +\xx')=
e^{i\phi^{\epsilon A}(\xx+\gamma,\xx' + \gamma')}\tilde{
\lambda} (\gamma -\gamma')\delta(\xx-\xx').$$
We see that $T_{\gamma\gamma'}$ is a multiplication operator:
$$L^2(\Omega)\ni f\mapsto [T_{\gamma\gamma'}f](\xx)=e^{i\phi^{\epsilon A}(\xx+\gamma,\xx + \gamma')}\tilde{
\lambda} (\gamma -\gamma') f(\xx)\in L^2(\Omega).$$
Consider the unitary operator $$U_\epsilon : l^2(\Gamma)\otimes L^2(\Omega)\mapsto l^2(\Gamma)\otimes L^2(\Omega),\quad 
[U_\epsilon \Psi]_\gamma(\xx):= e^{i\phi^{\epsilon A}(\gamma,\xx+\gamma)} \Psi_\gamma(\xx). $$
The operator 
$$T^\epsilon := U_\epsilon \;\mathfrak{Op}^{\epsilon A}(\lambda) \; U_\epsilon^*,\quad [T^\epsilon_{\gamma \gamma'}f](\xx)=e^{i\phi^{\epsilon A}(\gamma,\xx+\gamma)}e^{i\phi^{\epsilon A}(\xx+\gamma,\xx + \gamma')}e^{i\phi^{\epsilon A}(\xx+\gamma',\gamma')}\tilde{
\lambda} (\gamma -\gamma') f(\xx)$$
will have the same spectrum as $\mathfrak{Op}^{\epsilon A}(\lambda)$. Define the operator 
$$ \widetilde{T}^\epsilon=\{\widetilde{T}^\epsilon_{\gamma, \gamma'}\}_{\gamma \gamma'\in \Gamma},\quad \widetilde{T}^\epsilon_{\gamma \gamma'}=e^{i\phi^{\epsilon A}(\gamma,\gamma')}\tilde{
\lambda} (\gamma -\gamma')\bb1.$$
The following estimate
\begin{align*} 
 &|\phi^{\epsilon A}(\gamma,\xx+\gamma)+\phi^{\epsilon A}(\xx+\gamma,\xx + \gamma')+\phi^{\epsilon A}(\xx+\gamma',\gamma')-\phi^{\epsilon A}(\gamma,\gamma')|\\
 &\leq |\epsilon| C \left (|\la \gamma, \xx+\gamma, \xx + \gamma'\ra | +|\la \gamma, \xx+\gamma',\gamma'\ra | \right )
 \end{align*}
is a consequence of \eqref{iunie26} applied twice. We observe that since $\Omega$ is bounded, the areas of both triangles are bounded from above by $|\gamma-\gamma'|$. Using a Schur-Holmgren type bound, we obtain that:
$$||T^\epsilon -\widetilde{T}^\epsilon||\leq C |\epsilon| \sum_{\gamma\in \Gamma} |\tilde{
\lambda} (\gamma)|\; |\gamma|$$
which shows that the spectrum of $\mathfrak{Op}^{\epsilon A}(\lambda)$ is at an $|\epsilon|$-Hausdorff distance from the spectrum of $\widetilde{T}^\epsilon$. Hence the spectral edges of  $\mathfrak{Op}^{\epsilon A}(\lambda)$ have the same regularity as those of $\widetilde{T}^\epsilon$. The operator $\widetilde{T}^\epsilon$ is independent of the $\xx$ variable and we have:
$$ \widetilde{T}^\epsilon=\widetilde{t}^\epsilon \otimes \bb1,\quad \widetilde{t}^\epsilon=\{\widetilde{t}^\epsilon_{\gamma \gamma'}\}_{\gamma, \gamma'\in \Gamma}\in \mathbb{B}(l^2(\Gamma)),\quad \widetilde{t}^\epsilon_{\gamma \gamma'}=e^{i\phi^{\epsilon A}(\gamma,\gamma')}\tilde{
\lambda} (\gamma -\gamma').$$

Hence it is enough to study the spectral edges of the discrete operator $\widetilde{t}^\epsilon$ acting on $l^2(\Gamma)$, which is exactly of the form previously considered in \cite{Ne2} and \cite{C}. Although the Lipschitz behavior up to the logarithmic factor is essentially proved in \cite{Ne2}, let us show how one can modify the proof of our Theorem \ref{theorem1} in order to cover the discrete case. 

First of all, the space $\mathscr{C}^\alpha$ introduced in Definition \ref{definition1} will now consist of operators $t\in \mathbb{B}(l^2(\Gamma))$ for which:
$$||t||_\alpha=\max\left \{\underset{\gamma\in\Gamma}{\sup}\sum_{\gamma'}|t_{\gamma\gamma'}|\la \gamma-\gamma'\ra ^\alpha,\; \underset{\gamma'\in\Gamma}{\sup}\sum_{\gamma\in \Gamma}\left|t_{\gamma\gamma'}\right|\la \gamma-\gamma'\ra ^\alpha\,d\x\right \}<\infty.
 $$
 Now if everywhere in the proof of Theorem \ref{theorem1} we replace $\x$ and $\x'$ with $\gamma$ and $\gamma'$, the Lebesgue integration over $\R^d$ with respect to $d\x$ and $d\x'$ with sums over $\Gamma$, and $L^2(\R^d)$ with $l^2(\Gamma)$, everything remains true. Note the important fact that the integration with respect to $\y$ {\it must not} be replaced with a discrete sum. 
 
 We conclude that the spectral edges (and the norm) of $\widetilde{t}^\epsilon$ (hence $\mathfrak{Op}^{\epsilon A}(\lambda)$) obey the estimates announced in Corollary \ref{corolar1}, where the constants are proportional with the quantity 
 $$\sum_{\gamma \in \Gamma} \langle \gamma\rangle ^2 |\tilde{
\lambda} (\gamma)|.$$
\qed

\section*{Acknowledgements}
H.C. was supported by Grant 11-106598 of the Danish Council for Independent Research $|$ Natural Sciences, and a  Bitdefender Invited Professor Scholarship with IMAR, Bucharest. R.P. acknowledges the partial support of a grant of the Romanian National
Authority for Scientific Research, CNCS-UEFISCDI, project number
PN-II-ID-PCE-2011-3-0131 and the hospitality of the Aalborg University where
part of this work has been done. Both authors thank Gheorghe Nenciu for many illuminating discussions.


\begin{thebibliography}{00}

\bibitem[ABG]{ABG} Amrein, W.O., Boutet de Monvel, A., Georgescu, V.:
{\it $C_0$-Groups, Commutator Methods and Spectral Theory of N-Body
Hamiltonians}, Birkh\"auser, Basel, 1996.


\bibitem[AMP]{AMP} Athmouni, N., M\u{a}ntoiu, M., Purice, R.:  On the
continuity of spectra for families of magnetic pseudodifferential operators. Journal of Mathematical Physics {\bf 51}, 083517 (2010); doi:10.1063/1.3470118 (15 pages)

\bibitem[AS]{AS}Avron, J.E., Simon, B.: Stability of gaps for periodic potentials under variation of a magnetic
field. J. Phys. A: Math. Gen. {\bf 18}, 2199-2205 (1985)

\bibitem[AMS]{AMS}Avron, J., van Mouche, P.H.M., Simon, B.: On the measure of the spectrum for the almost
Mathieu operator. Commun. Math. Phys. {\bf 132}, 103-118, (1990). 
Erratum in Commun. Math.
Phys. {\bf 139}, 215 (1991)

\bibitem[B1]{Bell1}Bellissard, J.: Le papillon de Hofstadter. S{\'e}minaire Bourbaki {\bf 34}, 7-39 (1991-1992)

\bibitem[B2]{Bell2}Bellissard, J.: Lipshitz Continuity of Gap Boundaries
for Hofstadter-like Spectra. Commun. Math. Phys. {\bf 160}, 599-613 (1994)

\bibitem[BKS]{BCS}Bellissard, J., Kreft, C., Seiler, R.: Analysis of the spectrum of a particle on a triangular lattice with two magnetic fluxes by algebraic and numerical methods. J. of Phys. A {\bf 24}(10), 2329-2353 (1991)


\bibitem[BC]{BC} Briet, P., Cornean, H.D.: 
Locating the spectrum for magnetic Schr\"odinger and Dirac operators.  
{\it Comm. Partial Differential Equations}  {\bf 27}  no. 5-6, 1079--1101 
(2002)

\bibitem[BEY]{BEY}Choi, M.D., Elliott, G.A., Yui, N.: 
Gauss polynomials and the rotation algebra. Invent. Math.
{\bf 99}, 225-246 (1990)

\bibitem[C]{C} Cornean, H.D.:  On the Lipschitz Continuity of Spectral Bands of
Harper-Like and Magnetic Schr\"odinger Operators. Annales Henri Poincar\'{e}
{\bf 11}(5), 973--990 (2010).


\bibitem[CN]{CN} Cornean, H.D., Nenciu, G.: 
\newblock On eigenfunction decay for two dimensional magnetic Schr\"odinger operators.
\newblock Commun. Math. Phys. {\bf 192}, 671-685 (1998) 

 \bibitem[CP]{Co-Pu} Cornean, H.D., Purice, R.: 
 On the Regularity of the Hausdorff Distance Between Spectra of Perturbed Magnetic Hamiltonians. Operator Theory: Advances and Applications {\bf 224}, 55--66 (2012) 


\bibitem[DGR]{DGR} Dombrowski, N., Germinet, F., Raikov, G.D.:
Quantization of edge currents along magnetic barriers and magnetic guides. 
Ann. H. Poincar{\'e} {\bf 12}, 1169-1197 (2011) 






\bibitem[E]{Ell}Elliott, G.: Gaps in the spectrum of an almost periodic 
Schrodinger operator. C.R. Math. Rep.
Acad. Sci. Canada {\bf 4}, 255-259 (1982)




\bibitem[HR]{HR}Haagerup, U., R{\o }rdam, M.: 
Perturbations of the rotation $C^*$-algebras and of the
Heisenberg commutation relation. Duke Math. J. {\bf 77}, 627-656 (1995)


\bibitem[HKS]{HKS} Helffer, B., Kerdelhue, P., Sj{\"o}strand, J.:
  M{\'e}moires de la SMF, S{\'e}rie 2 {\bf 43}, 1-87 (1990)


\bibitem[HS1]{He-Sj1} Helffer, B., Sj{\"o}strand, J.: 
\newblock Equation de Schr\"odinger avec champ magn{\'e}tique et {\'e}quation de Harper.
\newblock  Springer Lecture Notes in Phys.  {\bf 345}, 118-197  (1989)

\bibitem[HS2]{He-Sj2} Helffer, B., Sj{\"o}strand, J.: Analyse
  semi-classique pour l'{\'e}quation de Harper. II. Bull. Soc. Math.
France {\bf 117}, Fasc. 4, Memoire 40 (1990)


\bibitem[Ho]{Ho-4} H\"{o}rmander, L.: {\it Analysis of Partial Differential Operators. Vol. IV}.




\bibitem[IMP1]{IMP1} Iftimie, V., M\u antoiu, M., Purice, R.: Magnetic
Pseudodifferential Operators. Publ. RIMS. {\bf 43}83), 585--623 (2007)

\bibitem[IMP2]{IMP2} Iftimie, V., M\u antoiu, M., Purice, R.:  Commutator Criteria
for magnetic pseudodifferential operators. Comm. Partial Diff. Equations. {\bf
35}, 1058--1094 (2010)

\bibitem[K]{Ko} Kotani, M.: 
Lipschitz continuity of the spectra of the magnetic transition
operators on a crystal lattice.  J. Geom. Phys.  {\bf 47} (2-3), 323--342 (2003)

\bibitem[LMR]{LMS} Lein, M., M\u antoiu, M., Richard, S.: Magnetic
  pseudodifferential operators with coefficients in $C^*$-algebras. Publ. RIMS Kyoto Univ. {\bf 46}, 755-788 (2010)

\bibitem[MP1]{MP1}M\u antoiu, M., Purice, R.: 
Strict deformation quantization for a particle in a magnetic field.  
 J. Math. Phys.  {\bf 46}  (5), 052105 (2005) 	


\bibitem[MP2]{MP2}M\u antoiu, M., Purice, R.: The magnetic Weyl calculus. 
 J. Math. Phys. {\bf 45} (4), 1394--1417 (2004)




\bibitem[MPR]{MPR} M\u antoiu, M., Purice, R., Richard, S.:  Spectral and
Propagation Results for Magnetic Schr\"odinger Operators; a
$C^*$-Algebraic Approach. J. Funct. Anal. {\bf 250}, 42--67  (2007)


\bibitem[N1]{Nenciu} Nenciu, G.: Stability of energy gaps under variation of the magnetic field. Lett. Math. Phys. {\bf 11}, 127-132 (1986)

\bibitem[N2]{Ne1} Nenciu, G.: On asymptotic perturbation theory for quantum
mechanics: Almost invariant subspaces and gauge invariant magnetic perturbation
theory. J. Math. Phys. {\bf 43}(3), 1273-1298 (2002) 

\bibitem[N3]{Ne2} Nenciu, G.: On the smoothness of gap boundaries for
generalized Harper operators. In: Advances in operator algebras and
mathematical physics, 173--182, Theta Ser. Adv. Math., 5, Theta, Bucharest, 2005.













\end{thebibliography}
\end{document}